\documentclass{amsart}
\usepackage{amssymb,amsmath,latexsym,amsthm}%or any usual package
\usepackage[english]{babel}

\usepackage{amsmath,amssymb,amsbsy,amsfonts,amsthm,latexsym,
                       amsopn,amstext,amsxtra,euscript,amscd}

\usepackage{array}
\usepackage{ifthen}
\usepackage{url}
\usepackage{hyperref}

\newtheorem{theorem}{Theorem}[section]
\newtheorem{lemma}{Lemma}[section]

\theoremstyle{definition}

\def\hh{{\rm h}\kern.4pt}

\begin{document}%51

\title{On Erd\H os' last equation}

\author{Florian~Luca}
\address{Mathematics Division, Stellenbosch University, Stellenbosch, South Africa
\newline \hspace*{4mm}
Max Planck Institute for Software Systems, Saarbr\"ucken, Germany
\newline \hspace*{4mm}
Department of Computer Science, Oxford University, UK
}
\email{fluca@sun.ac.za}

\author{Istv\'an~Pink}
\address{Institute of Mathematics, University of Debrecen,\newline
\indent P. O. Box 400, H-4002 Debrecen, Hungary}
\email{pinki@science.unideb.hu}

\author{Csaba~S\'andor}
\address{Department of Stochastics, Institute of Mathemetics, Budapest University of Technology and Economics, M\H{u}egyetem rkp. 3., H-1111 Budapest, Hungary; \newline \hspace*{4mm}
HUN-REN Alfr\'ed R\'enyi Institute of Mathematics, Re\'altanoda utca 13--15., H-1053 Budapest, Hungary
MTA--HUN-REN RI Lend\"ulet ``Momentum'' Arithmetic Combinatorics Research Group, Re\'altanoda utca 13--15., H-1053 Budapest, Hungary}
\email{sandor.csaba@ttk.bme.hu}

\thanks{Florian Luca was partially supported by the 2024 ERC Synergy Grant ``DynAMiCs. The second author is supported by the grants 150284 and ANN130909 of the Hungarian
National Foundation for Scientific Research. The research was supported by the University of
Debrecen Program for Scientific Publication.
The third author is supported by the NKFIH Grants No. K129335, KKP 144059 and the Lend\"ulet ''Momentum'' program of the Hungarian Academy of Sciences (MTA)}

\date{\today}

\pagenumbering{arabic}

\maketitle

\begin{abstract}
For a positive integer $n\ge 2$, we look at the number of solutions of the Diophantine equation
$$
n(x_1+\cdots+x_n)=x_1\cdots x_n.
$$
Denoting by $f(n)$ this number we show that $\liminf_{n\to\infty} \log f(n)/\log n=0$.
\end{abstract}

\section{Introduction}

On 19 August 1996, a month before his death, Paul Erd\H os formulated the following question in a letter to Richard Guy \cite{Guy}:
``Do you know anything about the equation
$$
n(x_1+x_2+\cdots+x_n)=x_1\cdots x_n
$$
$1\le x_1\le \cdots \le x_n$, where the $x_i$ are positive integers? Denote the number of solutions by $f(n)$; $f(n)>n^{\varepsilon}$ for some $\varepsilon>0$ seems to be true? \ldots".

This equation has been called Erd\H os' last equation. It has been studied in a couple of papers (see  \cite{LPS} and \cite{Shiu}). Yet, none of these papers addressed Erd\H os' claim.
In this paper, we show that the answer to Erd\H os' question is in the negative. We have the following result.

\begin{theorem}
We have
$$
\liminf_{n\to \infty} \frac{\log f(n)}{\log n}=0.
$$
\end{theorem}

It is not hard to show that
$$\limsup \frac{\log f(n)}{\log n}\le 2.
$$
We believe that
\begin{equation}
\label{eq:equal1}
\limsup_{n\to\infty} \frac{\log f(n)}{\log n}=1,
\end{equation}
and we will present an heuristic supporting this conjecture.
It is easy to show that
$$\limsup_{n\to\infty} \frac{\log f(n)}{\log n}\ge 1,
$$
and we will indicate a quick proof of this, so for \eqref{eq:equal1} one only needs to show that the limsup is at most $1$. Restricting to primes, we prove that
$$\liminf_{p\to\infty} \frac{\log f(p)}{\log p}=0,\qquad \limsup_{p\to\infty} \frac{\log f(p)}{\log p}\le 1\qquad {\text{\rm  and}}\qquad
\limsup_{p\to\infty} \frac{\log f(p)}{\log p}\ge \frac{50}{259}.
$$
In particular, writing for $n\ge 2$, $f(n)=n^{c(n)}$, we have that $\liminf_{n\to\infty} c_n=0$ and $\limsup_{n\to\infty} c_n\in [1,2]$.  Thus, the set of cluster points of $\{c_n\}_{n\ge 1}$ is in $[0,2]$. We leave the question of describing the cluster points for $\{c_n\}_{n\ge 1}$ as a challenge to the reader.

\section{Preliminary lemmas}

We rewrite the equation as
$$
x_n=\frac{n(x_1+\cdots+x_{n-1})}{x_1\cdots x_{n-1}-n}.
$$
Since $x_i$ are positive integers, it follows that $m:=x_1\cdots x_{n-1}-n$ is a positive integer. We will find some bounds for $m$.

\begin{lemma}
\label{lem:1}
We have $m<n(n-1)$ for $n\ge 7$.
\end{lemma}

\begin{proof}
We have
\begin{equation}
\label{eq:1}
x_1+\cdots+x_{n-1}\le (n-1)x_n.
\end{equation}
In particular, we have
$$
x_n=\frac{n(x_1+\cdots+x_{n-1})}{m}\le \frac{n(n-1)x_n}{m},
$$
which leads to $m\le n(n-1)$. The inequality is strict assuming that \eqref{eq:1} is strict, which is the case except when $x_1=\cdots=x_{n-1}=x_n$. In this last case, denoting the common value of the $x_i$'s by $x$, we get
$$
x=\frac{n(n-1)x}{x^{n-1}-n}\qquad {\text{\rm which leads to}}\qquad x^{n-1}=n+n(n-1)=n^2.
$$
Since $x\ge 2$, we get $n^2\ge 2^{n-1}$, which is false for $n\ge 7$. Hence, the inequality is strict for $n\ge 7$.
\end{proof}

We now show a different inequality for $m$.
\begin{lemma}
\label{lem:2}
Assume $n\ge 7$. Let $d:=\gcd(n,m)$. If $m\ge 2dn$, then $m < 8n\log n$.
\end{lemma}

\begin{proof}
Write $n:=n_1d,~m:=m_1d$. Assume $m_1\ge 2n$. Then
$$
x_n=\frac{n(x_1+\cdots+x_{n-1})}{m}=\frac{n_1(x_1+\cdots+x_{n-1})}{m_1}.
$$
Hence, $n_1\mid x_n$ and
\begin{equation}
\label{eq:imp1}
m_1\left(\frac{x_n}{n_1}\right)=x_1+\cdots+x_{n-1}.
\end{equation}
We also have $x_1\cdots x_{n-1}=m+n$. Let $k$ be the number of $x_i$'s which are larger than  $1$. Putting $j:=n-k-1$, we have $1=x_1=\cdots=x_j$ and
$2\le x_{j+1}\le \cdots\le x_{n-1}$. Then
\begin{equation}
\label{eq:imp2}
x_1+\cdots+x_{n-1}=n-k-1+x_{j+1}+\cdots+x_{n-1}
\end{equation}
is a multiple of $m_1$. This shows that
$$
kx_{n-1}\ge x_{j+1}+\cdots+x_{n-1}=(x_1+\cdots+x_{n-1})-(n-k-1)> m_1-n\ge \frac{m_1}{2},
$$
so
$$
x_{n-1}\ge\frac{m_1}{2k}.
$$
On the other hand, we have
$$
\frac{m}{n} x_n=\frac{m_1}{n_1} x_n=x_1+\cdots+x_{n-1}\le n-k-1+kx_{n-1}\le n-k-1+kx_n,
$$
which gives
\begin{equation}
\label{eq:31}
\left(\frac{m}{n}-k\right)x_n\le n-k-1.
\end{equation}
We have, by Lemma \ref{lem:1},
\begin{equation}
\label{eq:4}
k\le \frac{\log(m+n)}{\log 2}\le \frac{\log(n^2)}{\log 2}<4\log n.
\end{equation}
If $m/n<8\log n$, then we get what we want. Thus, we may assume on the contrary that $m/n\ge 8\log n>2k$, which together with \eqref{eq:31} shows that
$$
\frac{m_1}{2k}\le x_{n-1}\le x_n\le \frac{2n(n-k-1)}{m}< \frac{2nn_1}{m_1}.
$$
This gives
$$
m_1^2<4kn n_1<16 nn_1\log n.
$$
We thus get that
$$
m_1<4{\sqrt{n n_1 \log n}}=\frac{4n}{\sqrt{d}} {\sqrt{\log n}},
$$
which gives
\begin{equation}
\label{eq:3}
m< 4n{\sqrt{d}} {\sqrt{\log n}}.
\end{equation}
However, we assumed that $m\ge 2dn$. We thus get
$$
2dn<4n{\sqrt{d}} {\sqrt{\log n}}\qquad {\text{\rm which gives}}\qquad {\sqrt{d}}< 2{\sqrt{\log n}},
$$
which together with \eqref{eq:3} gives $m< 8n\log n$, which is what we wanted.
\end{proof}

The combination of Lemma \ref{lem:1} and \ref{lem:2} gives good information about the shape of the solutions to Erd\H os' equation when $n=p$ is prime.
Indeed, in this case, $d\in \{1,p\}$. If $d=1$, then Lemma \ref{lem:2} shows that $m<8p\log p$, while if $d=p$, then Lemma \ref{lem:1} shows that
$m=pm_1$, where $m_1<p-1$ for $p\ge 7$.

\section{The proof}

For a positive integer $n$ we write $F(n)$ for the number of unordered factorizations of $n$ in factors $>1$ with the convention that $F(1)=1$. That is, for $n>1$,
$$
F(n):=\#\{2\le n_1\le \cdots\le n_k: n_1\cdots n_k=n\}.
$$
If we want to emphasize the number of factors in the factorization, we write $F_k(n)$ for the number of factorizations of $n$ with exactly $k$ factors. Note that
$$
F(n)=\sum_{k\ge 1} F_k(n).
$$
In the above sum, $F_k(n)$ is zero once $k>\Omega(n)$. We now return to Erd\H os' equation for $n=p$. With the notation from the proof of Lemma \ref{lem:2}, any solution is of the form
\begin{equation}
\label{eq:z}
z=\frac{p(p-k-1+y_1+\cdots+y_k)}{m},
\end{equation}
where $2\le y_1\le \cdots\le y_k$, $y_1\cdots y_k=p+m$. Here, $y_1,\ldots,y_k$ are the integers denoted by $x_{j+1}\ldots, x_p$ in Lemma \ref{lem:2} and
$z$ is the integer denoted $x_p$. We write
$$f(p)=:f_1(p)+f_2(p),$$
where $f_1(p)$ accounts for those solutions for which $\gcd(m,p)=1$ and
$f_2(p)$ accounts for those solutions for which $\gcd(m,p)=p$. For a large real number $X$ we will upper bound
$$
S_i(X):=\sum_{X/2\le p\le X} f_i(p)\qquad {\text{\rm for}}\qquad i=1,2.
$$
We will show that
$$
S_i(X)=X^{1+o(1)}\qquad {\text{\rm as}}\qquad X\to\infty.
$$
This shows that
$$
\sum_{X/2\le p\le X} f(p)=S_1(X)+S_2(X)=X^{1+o(1)}\qquad {\text{\rm as}}\qquad X\to\infty.
$$
Since the interval $[X/2,X]$ contains $\pi(X)-\pi(X/2)=X^{1+o(1)}$ primes, it follows that
$$
\min_{X/2\le p\le X} f(p)\le X^{o(1)}=p^{o(1)}\qquad {\text{\rm as}}\qquad X\to\infty,
$$
which proves our theorem. It also shows that $\limsup_{p\to\infty} \log f(p)/\log p\le 1$.
Let's get to work.

\subsection{Estimating $S_1(X)$}

In this case, if $p\in [X/2,X]$ and $X\ge 14$, we then have $m\le 8p\log p\le 8X\log X$ by Lemma \ref{lem:2}. Thus, $p+m\le p+8p\log p\le 9X\log X$
for $X\ge 14$. Fix $k$, which is the number of terms larger than $1$ in  the factorization $y_1\cdots y_k=p+m$. Certainly, $k\le \Omega(p+m)=O(\log X)$, where $\Omega(.)$ counts the total number of prime factors with multiplicity.
The number of such factorizations is $F_k(p+m)$.  In equation \eqref{eq:z}, we have that $p$ divides $z$, so $z=pz_1$ and
\begin{equation}
\label{eq:z1}
z_1=\frac{p-k-1+y_1+\cdots+y_k}{m}=\frac{(p+m)-k-1+y_1+\cdots+y_k}{m}-1,
\end{equation}
so $m$ is a divisor of $y_1\cdots y_k-k-1+y_1+\cdots+y_k$. Note that this last number is positive as $k\ge 1$ since $y_1\cdots y_k=p+m>1$.
This suggests that any pair $(p,x_1,\ldots,x_p)$ of solutions of the type participating in $S_1(X)$ of Erd\H os' equations arises from the following parameters:
\begin{itemize}
\item[(i)] A positive integer $k$ of size $O(\log X)$.
\item[(ii)] A positive integer $\ell(:=p+m)$ which is at most $O(X\log X)$.
\item[(iii)] An unordered factorization $y_1\cdots y_k=\ell$ with $k$ factors larger than $1$ of $\ell$.
\item[(iv)] A proper divisor $m$ of the positive integer $y_1\cdots y_k-k-1+y_1+\cdots+y_k$.
\end{itemize}
The condition in (iv) that $m$ is a ``proper" divisor is meant to insure that the number $z_1$ arising in \eqref{eq:z1} is not zero.
Indeed, once (i)--(iv) are known, then $p=\ell-m$, and the solution is $x_1=\cdots=x_j=1$, $(x_{j+1},\ldots,x_{p-1})=(y_1,\ldots,y_{k})$ and
$x_p=pz_1$, where $z_1$ is shown in \eqref{eq:z1}. Note that this might not give rise to a genuine solution since the $p$ arising from this process is not necessarily prime.
This shows that
$$
S_1(X)\ll \sum_{\ell\le 9X\log X} \sum_{k\le \Omega(\ell)} \sum_{\substack{2\le y_1\le \cdots\le y_k\\ y_1\cdots y_k=\ell}} \tau(y_1\cdots y_k-k-1+y_1+\cdots+y_k),
$$
where $\tau(.)$ counts the number of positive divisors.
In the inner sum, all numbers whose number of divisors is being calculated are at most $\ell+ky_k\le \ell+k\ell/2^{k-1}\le 2\ell\le 18X\log X$ so
each of these numbers is $X^{o(1)}$ as $X\to\infty$ and the number of them is $F_k(\ell)$. Thus,
\begin{equation}
\label{eq:S1}
S_1(X)\le X^{o(1)} \sum_{\ell\le 9X\log X} F(\ell).
\end{equation}
By a well-known result of Szekeres and Tur\'{a}n \cite{ST}, we have
\begin{equation}
\label{eq:ST}
\frac{1}{t}\sum_{n\le t} F(n)=(1+o(1)) \frac{\exp(2{\sqrt{\log t}})}{2{\sqrt{\pi}} (\log t)^{3/4}}\qquad {\text{\rm as}}\qquad t\to\infty.
\end{equation}
Inserting \eqref{eq:ST} into \eqref{eq:S1} yields
$$
S_1(X)\ll  X^{1+o(1)}(\log X) \frac{\exp(2{\sqrt{\log(9X\log X)}})}{2{\sqrt{\pi}}(\log(9X\log X))^{3/4}}=X^{1+o(1)}\quad {\text{\rm as}}\qquad X\to\infty,
$$
as desired.

\subsection{Estimating $S_2(X)$}
We prove that $f_2(p)=3$ for $p\ge 3$. It is enough to verify the following lemma.

\begin{lemma}
\label{lem:3}
Let $p$ be an odd prime number. The solutions $(x_1,\ldots,x_p)$ of the equation
$$
p(x_1+x_2+\dots +x_p)=x_1x_2\dots x_p\qquad {\text{\rm with}}\qquad p\mid x_1x_2\dots x_{p-1}
$$
are $(1,1,\dots ,1,2p,3p-2)$, $(1,1,\dots ,1,2,p,2p-1)$ and $(1,1,\dots ,1,3,p,p)$.
\end{lemma}
\begin{proof}
First we show that $p\nmid x_1x_2\dots x_{p-2}$ for $p\ge 5$. Let us suppose on the contrary that $p\mid x_1x_2\dots x_{p-2}$. Then $p\le x_{p-2}\le x_{p-1}\le x_p$ and $x_1+x_2+\dots +x_p<px_p$. We then have
$$
p^2x_p\le x_1x_2\dots x_{p-2}x_{p-1}x_p=p(x_1+x_2+\dots +x_p)<p^2x_p,
$$
a contradiction. So, we may assume that $x_{p-1}=kp$ for some positive integer $k$. We verify that $x_1x_2\cdots x_{p-2}\le 3$. Assume that $A=x_1x_2\cdots x_{p-2}\ge 4$. Since $$p(x_1+x_2+\dots +x_{p-2}+kp+x_p)=x_1x_2\cdots x_{p-2}kpx_p$$ then
$$
x_1+x_2+\cdots+ x_{p-2}+kp=(kA-1)x_p.
$$
Since
$$x_1+x_2+\dots +x_{p-2}\le p-3+x_1x_2\cdots x_{p-2}=p-3+A,
$$
we get
$$(kA-1)x_p\le p-3+A+kp,\qquad {\text{\rm that is}}\qquad x_p\le \frac{p-3+A+kp}{kA-1}.
$$
We get a contradiction if we show that
$$\frac{p-3+A+kp}{kA-1}<p.
$$
This is equivalent to $A-3<p(k(A-1)-2).$ This last inequality is equivalent to
$$
A-3<pk(A-1)-2p\quad {\text{\rm or}}\quad 2(p-1)<(A-1)(pk-1).
$$
The last inequality is equivalent to
$$
A-1> \frac{2(p-1)}{pk-1},
$$
which holds for $A\ge 4,~p\ge 2,~k\ge 1$. This contradiction shows that $A\in \{1,2,3\}$.

For $A=1$, $x_1=x_2=\dots =x_{p-2}=1$, and $p(p-2+x_{p-1}+x_p)=x_{p-1}x_p$ that is $(x_p-p)(x_{p-1}-p)=p(2p-2)$, which yields the solution
$$
(1,1,\dots ,1,2p,3p-2).
$$
For $A=2$, $x_1=x_2=\dots =x_{p-3}=1$, $x_{p-2}=2$ and $p(p-1+x_{p-1}+x_p)=2x_{p-1}x_p$ that is $(2x_p-p)(2x_{p-1}-p)=p(3p-2)$, which yields the solution
$$
(1,1,\dots ,1,2,p,2p-1).
$$
For $A=3$, $x_1=x_2=\dots =x_{p-3}=1$, $x_{p-2}=3$ and $p(p+x_{p-1}+x_p)=3x_{p-1}x_p$ that is $(3x_p-p)(3x_{p-1}-p)=4p^2$, which yields the solution
$$(1,1,\dots ,1,3,p,p).
$$
\end{proof}

By Lemma \ref{lem:3}
$$
S_2(X)=\sum_{X/2\le p\le X} f_2(p)=\sum_{X/2\le p\le X} 3\le 3X\le X^{1+o(1)} \quad {\text{\rm as}}\quad X\to\infty,
$$
as desired.

\section{Comments}

In this section, we justify some of the claims made in the introduction about the $\limsup$ and $\liminf$ of $\log f(n)/\log n$. To get an upper bound, we may proceed in the same was as in the proof of the theorem. Namely, let $k\ge 1$ and $y_1\cdots y_k=n+m$ be a factorization of $\ell=n+m$ with $k$ factors larger than $1$. Putting $n\in [X/2,X]$ for $X\ge 14$, we have $n+m\le X^2$
so $k\le \Omega(n+m)\le 4\log X$. Next
$$
z=\frac{n(n-k-1+\sum_{i=1}^k y_i}{m}=n\left(\frac{(n+m)-k-1+\sum_{i=1}^k y_i}{m}-1\right).
$$
We have that $m=dm_1$, where $m_1$ is a proper divisor of $y_1\cdots y_k-k-1+y_1+\cdots+y_k$, and $d$ is a divisor of both $n$ and $m$; hence, of
$\ell=n+m=y_1\cdots y_k$. Thus,
$$
\sum_{X/2\le n\le X} f(n)\ll \sum_{\ell\le X^2} \sum_{k\le 4\log X} \sum_{\substack{2\le y_1\le \cdots\le y_k\\ y_1\cdots y_k=\ell} } \tau(\ell)\tau(\ell-k-1+y_1+\cdots+y_k).
$$
Each of the terms of the inner summation is of size $X^{o(1)}$ and the number of such terms is $F_k(\ell)$. Thus,
$$
\sum_{X/2\le n\le X} f(n)\ll \sum_{\ell\le X^2} \sum_{k\le 4\log X} F_k(\ell)\ll X^{2+o(1)},
$$
This shows that $f(n)\le n^{2+o(1)}$ as $n\to\infty$.  This is perhaps not the right order since the range for $m$ seems to be too high. As we said, we believe $f(n)\le n^{1+o(1)}$ holds
as $n\to\infty$. Here is an heuristic. As we saw, any solution of Erd\H os' equation is intimately related to unordered factorizations of the number $\ell:=n+m$. After such a factorization $y_1\cdots y_k=n+m$ with $k$ factors larger than $1$ is chosen, the quantity
$$
\frac{n(n-k-1+\sum_{i=1}^k y_i)}{m}
$$
needs to be an integer. Assuming that the probability that $n(n-k-1+\sum_{i=1}^k y_i)$ to be a multiple of $m$ to be $1/m$, it would seem that an upper bound on $f(n)$ should be
$$
\ll \sum_{m<n(n-1)} \frac{F(n+m)}{m}.
$$
We split this sum in two ranges, namely $m\le n$ and $m>n$. In the lower range, we ignore the denominator $m$ and say that the above sum is bounded by
$$
\sum_{m\le n} F(n+m)\le \sum_{\ell\le 2n} F(\ell)\le n^{1+o(1)}\qquad {\text{\rm as}}\qquad n\to\infty.
$$
In the larger range $m>n$, we have that $m>(n+m)/2$, so the sum in this range is bounded by
$$
\ll \sum_{n<m\le n^2} \frac{F(n+m)}{n+m}\ll \sum_{\ell\le n^2} \frac{F(\ell)}{\ell}\ll \exp(2{\sqrt{2\log \ell}})(\log \ell)^{1/4},
$$
where the last estimate follows by the Abel summation formula applied to \eqref{eq:ST}. The right--hand side above is $n^{o(1)}$ as $n\to\infty$. This seems to justify the claim that perhaps $f(n)\le n^{1+o(1)}$ as
$n$ tends to infinity. To see why $f(n)\ge n^{1+o(1)}$ holds for infinitely many $n$, we just take $m=1$. Then
$$
z:=n(n-k-1+\sum_{i=1}^k y_i)
$$
and we see right away that indeed $z:=x_n$ is an integer and larger than any $x_i$ for $i=1,\ldots,n-1$. Thus, $f(n)\ge F(n+1)$. Now choose
$n+1$ to be a highly--factorable integer in the terminology of \cite{CEP}, namely a positive integer $m$ such that $F(\ell)<F(m)$ for all $\ell<m$.
The paper \cite{CEP} shows that for such numbers, $F(m)=m L(m)^{-1+o(1)}$ as $m\to\infty$, where
$$
L(m):=\exp\left(\frac{\log m \log\log\log m}{\log\log m}\right)=m^{o(1)}\qquad {\text{\rm as}}\qquad m\to\infty.
$$
In particular, $f(n)\ge n^{1+o(1)}$ as $n$ tends to infinity in such a way that $n+1$ is highly--factorable. Finally, let us look at $\limsup_{p\to\infty} \log f(p)/\log p$.
If there were infinitely many primes $p$ such that $p+1$ is highly--factorable, then we would have that $f(p)\ge p^{1+o(1)}$ holds as $p$ tends to infinity through such numbers.
But this is not known. However, we can use known results on Linnik's problem which asks for the first (smallest) prime in an arithmetic progression.
For this, let $X$ be large and choose a highly--factorable number $m$. By Linnik's result, there is a prime $p$ of size $O(m^L)$ such that $m\mid p+1$ where $L$ is Linnik's constant. This shows that
$f(p)\ge F(p+1)\ge F(m)=m^{1/L+o(1)}$. The best (smallest) exponent $L$ is due to Xylouris \cite{Xy}  and is $<5$; hence, $1/L>1/5=0.2$.  To justify the details, we first
take $y:=C_1X^{1/L}$, where $C_1$ is an appropriate positive constant to guarantee that if $m\le y$ is any positive integer then there is a prime $p\le X$ such that
$p\equiv -1\pmod m$, and next we use the fact that there is a highly--factorable number in $[y/2,y]$ (see (2) on page 25 in \cite{CEP}). This justifies that
$f(p)\ge p^{0.2}$ holds as $p$ tends to infinity through some infinite subset of all the primes.


\begin{thebibliography}{9999}

\bibitem{CEP} E. R. Canfield, P. Erd\H os, and C. Pomerance, ``On a problem of Oppenheim concerning ``Factorization Numerorum"", {\it J. Number Theory} {\bf 17} (1983), 1--28.

\bibitem{Guy} R. K. Guy, {\em Unsolved Problems in Number Theory\/}, 3rd ed, New York: Springer--Verlag, 2004.

\bibitem{LPS} F. Luca, I. Pink and C. S\'andor, ``On the largest value of the solutions of Erd\H os' last equation", {\it Internat. J. Number Theory\/}, to appear.

\bibitem{ST} G. Szekeres and P. Tur\'{a}n, ``\"Uber das zweite Hauptproblem der "Factorisatio Numerorum"", Acta Litt. Sci. Szeged {\bf 6} (1933), 143--154.

\bibitem{Shiu} P. Shiu, ``On Erd\H os's last equation", {\it Amer. Math. Monthly\/} {\bf 126} (2019), no. 9, 802--808.

\bibitem{Xy} T. Xylouris, ``Linnik's constant is less than 5", {\it Chebyshevski{\u{\i}} Sb \/} {\bf 19} (2018), 80--94.

\end{thebibliography}
\end{document}